\documentclass[11pt,reqno]{amsart}
\usepackage[margin=1in]{geometry}
\usepackage[english]{babel}
\usepackage{amsmath}
\usepackage{amssymb}
\usepackage{amsthm}
\usepackage{latexsym}
\usepackage{mathtools}
\usepackage{amsfonts}

\usepackage{tikz-cd}

\usepackage{mathrsfs}
\usepackage{textcomp}
\usepackage{diagbox}
\usepackage[T1]{fontenc}
\usepackage{graphicx}
\usepackage{setspace}
\usepackage{nicefrac}
\usepackage{enumerate}
\usepackage{wasysym}
\usepackage{upgreek}
\usepackage{booktabs}
\usepackage{etoolbox}
\usepackage{wrapfig}
\usepackage{floatflt}
\usepackage{parskip}
\usepackage{lscape}
\usepackage{subcaption}

\usepackage[none]{hyphenat}
\usetikzlibrary{tikzmark,decorations.pathreplacing}
\usetikzlibrary{arrows.meta, positioning, shapes.geometric}
\usepackage{xcolor}
\definecolor{Maroon}{RGB}{128,0,0}
\usepackage{tikz}
\usepackage{colortbl}
\usepackage[hidelinks]{hyperref}
\usepackage{tikz-cd}
\usetikzlibrary{patterns,external}
\usetikzlibrary{calc}
\usepackage[numbers,sort]{natbib}
\usepackage{paralist}
\usepackage{mathtools}
\usetikzlibrary{automata}
\usepackage{tikz}
\usetikzlibrary{arrows,matrix}
\usetikzlibrary{positioning}
\usetikzlibrary{decorations}
\usepackage{cleveref}
\usepackage[all,cmtip]{xy}
\usepackage[labelfont=bf, font={small}]{caption}[2005/07/16]

\usepackage{mathrsfs}
\usepackage{textcomp}
\usepackage[T1]{fontenc}
\usepackage{graphicx}
\usepackage{setspace}
\usepackage{stackengine}
\usepackage{nicefrac}
\usepackage{enumerate}
\usepackage{wasysym}
\usepackage{upgreek}
\usepackage{hyperref}
\usepackage{booktabs}
\usepackage{wrapfig}
\usepackage{etoolbox}
\usepackage{subcaption}

\usetikzlibrary{decorations.pathmorphing,shapes}

\theoremstyle{definition}

\usepackage{aliascnt}
\usepackage{hyperref}

\newaliascnt{example}{theorem}
\newtheorem{example}[example]{Example}
\aliascntresetthe{example}

\usepackage[linesnumbered,commentsnumbered,ruled,vlined]{algorithm2e}

\definecolor{juliapurple}{RGB}{95,60,128}
\definecolor{juliagreen}{RGB}{56,152,38}
\definecolor{codegray}{RGB}{247,247,247}
\definecolor{darkteal}{RGB}{0,70,70}
\definecolor{juliaprompt}{RGB}{0,140,50}
\usepackage{listings}
\lstdefinelanguage{Julia}{
  morekeywords=[1]{
    using,import,end,return,if,else,elseif,for,while,in,do,
    struct,mutable,where,begin,let,const,nothing,true,false, NTuple, Float64, AbstractVector, Bool, push!, 
  },
  morekeywords=[2]{
    TriangulationCache, visualize, refine!, Sierpinski, Barycenter, Random, save
  },
  morekeywords=[3]{
    resolution, xlims, ylims, edges, buttons, strategy, is_complete, min_refinement_area
  },
  sensitive=true,
  morecomment=[l]{\#},
  morestring=[b]",
}
\newcounter{sarrow}

\makeatletter
\let\origsection\section
\renewcommand\section{\@ifstar{\starsection}{\nostarsection}}

\newcommand\nostarsection[1]
{\sectionprelude\origsection{#1}\sectionpostlude}

\newcommand\starsection[1]
{\sectionprelude\origsection*{#1}\sectionpostlude}

\newcommand\sectionprelude{%
  \vspace{-0.6em}
}

\newcommand\sectionpostlude{%
  \vspace{0em}
}
\makeatother

\title{Adaptive visualization for Enumerative Problems}
\author{Taylor Brysiewicz}
\email{tbrysiew@uwo.ca}
\address{Department of Mathematics,
University of Western Ontario,
London, Ontario, Canada}

\author{Noah Vale}
\email{nvale3@uwo.ca}
\address{Department of Mathematics,
University of Western Ontario,
London, Ontario, Canada}

\newcommand{\x}{\mathbf{x}}
\newcommand{\p}{\mathbf{p}}



\begin{document}
\begin{abstract}
We present an adaptive visualization framework and its implementation in the \texttt{julia} package \texttt{AdaptiveVisualization.jl}. Our algorithm visualizes  functions using an adaptive sampling procedure that concentrates the evaluation budget at interesting regions of the domain. We illustrate our package on our main application: visualizing parameter spaces of enumerative problems.   
\end{abstract}
\maketitle

\vspace{-20pt}
\section{Introduction}
The main goal of the \texttt{julia} \cite{julia} package \texttt{AdaptiveVisualization.jl} \cite{AdaptiveVisualizationjl} is to visualize a function
\[
f:\mathbb{R}^2\to\mathbb{R}
\]
given an evaluation oracle for $f$. The cost of this task is given by the number of oracle calls. Rather than sampling a compact domain $D\subset\mathbb{R}^2$ uniformly, our algorithm adaptively  samples in regions where the behavior of $f$ varies significantly. Therefore, the computational effort is focused on regions where the function has interesting behavior. 

\begin{figure}[!htpb]
\includegraphics[scale=0.5]{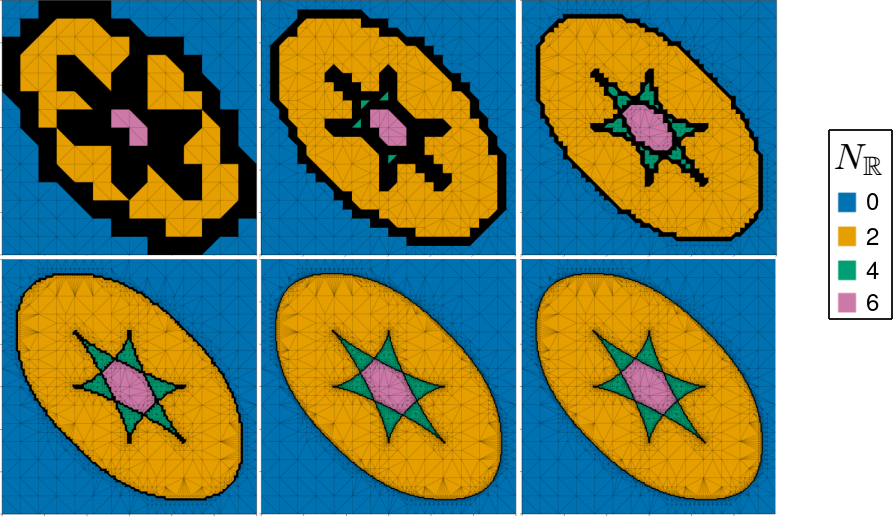} \vspace{-0.15in}
\caption{Successive refinements of visualizations of $N_{\mathbb{R}}$ for a  parametrized polynomial system with $6$ complex solutions: the equilibria of the $n=3$ Kuramoto model \cite{Kuramoto,Kuramoto1975}.}
\label{fig:kuramoto}
\end{figure}

Our framework is especially effective for functions associated with solution sets of parametrized polynomial systems. In this setting, an oracle is available by solving via homotopy continuation. For example, fix a polynomial system with parameters $p_1$ and $p_2$ which, for generic parameter values, admits finitely many solutions. This is known as an \textit{enumerative problem}. Define the function
\[
N_{\mathbb{R}}:\mathbb{R}_{p_1,p_2}^2 \to \mathbb{N} \hookrightarrow \mathbb{R}
\]
which counts the number of real solutions. It is locally constant off the \textit{real discriminant} $\Delta_{\mathbb{R}}$. Visualizing $N_{\mathbb{R}}$ shows the complement of $\Delta_{\mathbb{R}}$ directly, and a representation of $\Delta_{\mathbb{R}}$ emerges indirectly as the negative-space.  \Cref{fig:kuramoto} shows several refinements of a visualization of such a function. 

Our adaptive procedure involves triangulating a finite set of samples in a region $D$, and subsequently determining which triangles are \textit{complete} versus \textit{incomplete}.  The incomplete triangles are refined by inserting more samples and the process iterates as-decided by the user. At any point, the triangulation may be visualized using \texttt{GLMakie.jl}, the OpenGL backend of \texttt{Makie.jl} \cite{Makiejl}.

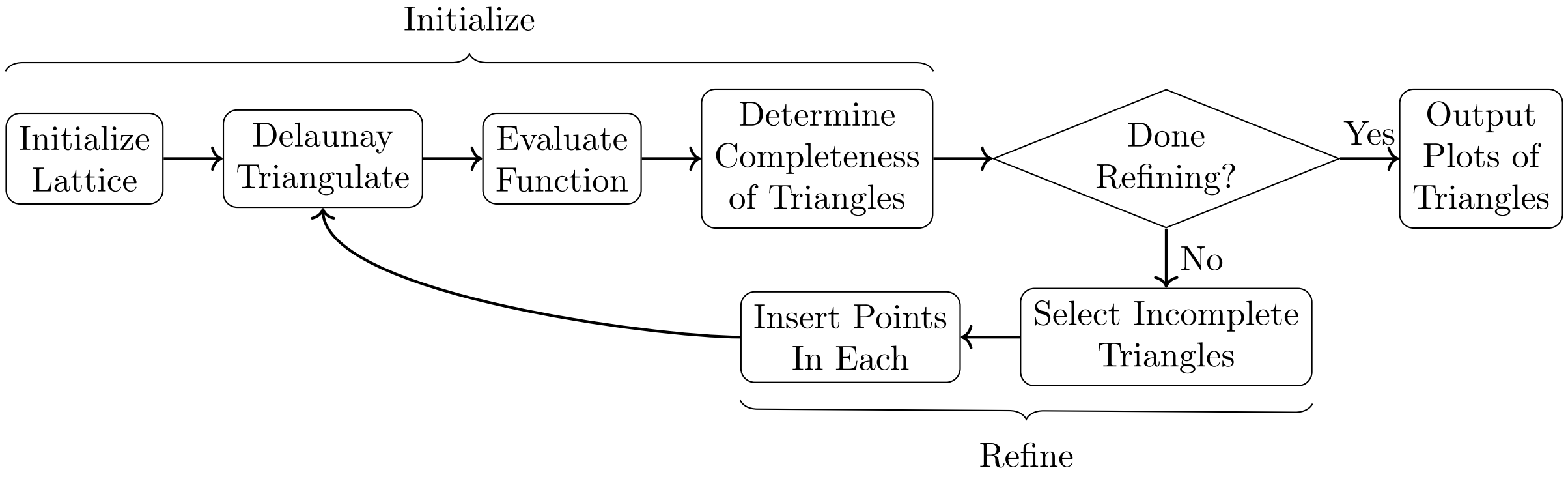
\begin{figure}[!htpb]
\centering
\begin{tikzpicture}[
    node distance=0.6cm,
    every node/.style={font=\small},
    process/.style={
        rectangle,
        draw,
        rounded corners,
        align=center,
        minimum width=0.8cm,
        minimum height=0.2cm
    },
    decision/.style={
        diamond,
        draw,
        align=center,
        aspect=2.5,
        inner sep=1pt
    },
    arrow/.style={->, thick}
]

% Main pipeline
\node (lattice) [process] {
    Initialize\\
    Lattice
};

\node (triangulate) [process, right=of lattice] {
    Delaunay\\
    Triangulate
};

\node (eval) [process, right=of triangulate] {
    Evaluate\\
    Function
};

\node (complete) [process, right=of eval] {
    Determine\\
    Completeness\\
    of Triangles
};

\node (done) [decision, right=of complete] {
    Done\\
    Refining?
};

\node (output) [process, right=of done] {
    Output\\
    Plots of\\
    Triangles
};

\node (foreach) [process, below=.6cm of done] {
    Select Incomplete \\
    Triangles
};

\node (insert) [process, left=of foreach] {
    Insert Points \\
    In Each
};

\coordinate (initialize-left)
    at (lattice.west |- complete.north);
\coordinate (initialize-right)
    at (complete.east |- complete.north);

\draw[
    decorate,
    decoration={brace,amplitude=5pt}
]
    ([yshift=5pt]initialize-left)
    --
    ([yshift=5pt]initialize-right)
    node[midway,above=8pt] {Initialize};

\draw[
    decorate,
    decoration={brace,mirror,amplitude=5pt}
]
    ([yshift=-5pt]insert.south west)
    --
    ([yshift=-5pt]foreach.south east)
    node[midway,below=8pt] {Refine};

% Main arrows
\draw [arrow] (lattice) -- (triangulate);
\draw [arrow] (triangulate) -- (eval);
\draw [arrow] (eval) -- (complete);
\draw [arrow] (complete) -- (done);
\draw [arrow]
    (done.east) -- node[above]{Yes} (output.west);

% Refinement branch
\draw [arrow]
    (done.south) -- node[right]{No} (foreach.north);

\draw [arrow]
    (foreach.west) -- (insert.east);

\draw [arrow]
    (insert.west)
    to[out=180,in=-90,looseness=0.5]
    (triangulate.south);

\end{tikzpicture}
\caption{Flowchart of the main steps of our adaptive visualization procedure.}
\label{fig:pipeline}
\end{figure}

Our software extends beyond real-valued functions on $\mathbb{R}^2$. Given an oracle for a function $f$ on $\mathbb{R}^n$, the software automatically restricts $f$ to a random $2$-plane $\Lambda$ and visualizes $f|_{\Lambda}$. It can distinguish discrete-valued functions like $N_{\mathbb{R}}$ from locally continuous functions like the \textit{Dietmaier function} \cite{Dietmaier}:
\begin{align*}
\mathcal D: \mathbb{R}^2 &\to \mathbb{R}\\
(p_1,p_2) &\mapsto\min_{\substack{\mathbf{x}\text{ a nonreal}\\\text{solution}}} (||\textrm{imag}(\textbf{x})||).
\end{align*}
Our software supports the visualization of functions which take categorical values, that is, discrete values which are not numerical.  For example, one may seek to simultaneously count the number of real and positive solutions to a parametrized system as a pair of integers.

In subsequent sections, we outline the pipeline shown in the flowchart of \Cref{fig:pipeline}. \Cref{sec:initialization} describes the initialization steps, and \Cref{sec:refinement}, the refinement steps. \Cref{sec:functionality} describes the user-interface. Finally, \Cref{sec:enumerative_problems} focuses on the \texttt{HomotopyContinuation.jl} extension for enumerative problems. 

\section{Initialization} 
\label{sec:initialization}

\setcounter{subsection}{-1}
We describe the procedure for initializing our main datatype: a \texttt{TriangulationCache}. To define it, we introduce two inputs to the initialization process: an evaluation oracle $\mathcal O$ for a function 
\begin{equation}
\label{eq:function}
f: \mathbb{R}^2 \to \mathbb{R}
\end{equation}
and a rectangular window 
\[
D = [x_{\textrm{min}},x_{\textrm{max}}]\times [y_{\textrm{min}},y_{\textrm{max}}].
\]
The next subsection describes how our software prepares more general functions into the form \eqref{eq:function}.

\begin{figure}[t]
$
\underset{\substack{\textbf{Initialize}\\\text{lattice}}}{
  \vcenter{\hbox{\includegraphics[scale=0.55]{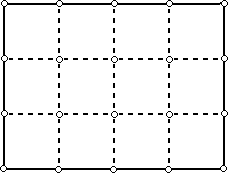}}}
}
\;\longrightarrow\;
\underset{\substack{\textbf{Triangulate}\\\text{points}}}{
  \vcenter{\hbox{\includegraphics[scale=0.55]{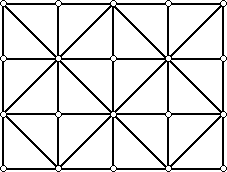}}}
}
\;\longrightarrow\;
\underset{\substack{\textbf{Evaluate}\\\text{oracle}}}{
  \vcenter{\hbox{\includegraphics[scale=0.55]{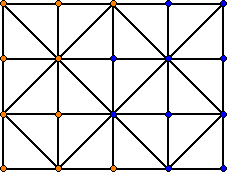}}}
}
\;\longrightarrow\;
\underset{\substack{\textbf{Completeness-check}\\\text{triangles}}}{
  \vcenter{\hbox{\includegraphics[scale=0.55]{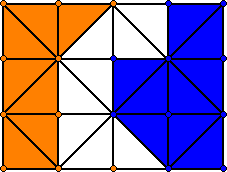}}}
}
$
\caption{Illustration of the four initialization steps.}
\end{figure}

\subsection{Pre-processing}
\label{sub:preprocess}
Given an evaluation oracle $\mathcal O$ for a function $f\vcentcolon \mathbb{R}^n \to \mathbb{R}$ with $n>2$, our software replaces $\mathcal O$ with the composition $\mathcal O \circ \Lambda$ where $\Lambda:\mathbb{R}_{x,y}^2 \to \mathbb{R}^n$ is a parametrization of a fixed random $2$-plane in $\mathbb{R}^n$. We sample a single random evaluation, and if it returns a floating point type, we declare that $f$ is (locally) \texttt{continuous}, and if it returns an integer, we declare it is \texttt{discrete}. More generally, if $f$ evaluates to a non-numerical value, we also declare that it is \texttt{discrete}. For the sake of  uniformity, discrete non-numerical values are interpreted as being embedded within $\mathbb{Z} \hookrightarrow \mathbb{R}$.

\subsection{Initial point mesh}
After pre-processing, the first step in initialization involves producing a lattice $S_0$ of points in the rectangle $D$ which preserves its aspect ratio and whose cells are approximately square. Given a resolution $N_0$, we choose the number of $x$ and $y$ values to be
\[
n_x =\left\lfloor \sqrt{N_0 \frac{x_{\textrm{max}}-x_{\textrm{min}}}{y_{\textrm{max}}-y_{\textrm{min}}}} \right\rfloor \quad \text{and} \quad n_y =\left\lfloor \sqrt{N_0 \frac{y_{\textrm{max}}-y_{\textrm{min}}}{x_{\textrm{max}}-x_{\textrm{min}}}}  \right\rfloor
\]
so that $n_xn_y = |S_0|  \approx N_0$ and $n_x/n_y \approx  \frac{x_{\textrm{max}}-x_{\textrm{min}}}{y_{\textrm{max}}-y_{\textrm{min}}}$.

\subsection{Delaunay Triangulation} 
\begin{figure}[!htpb]
\includegraphics[scale=0.25]{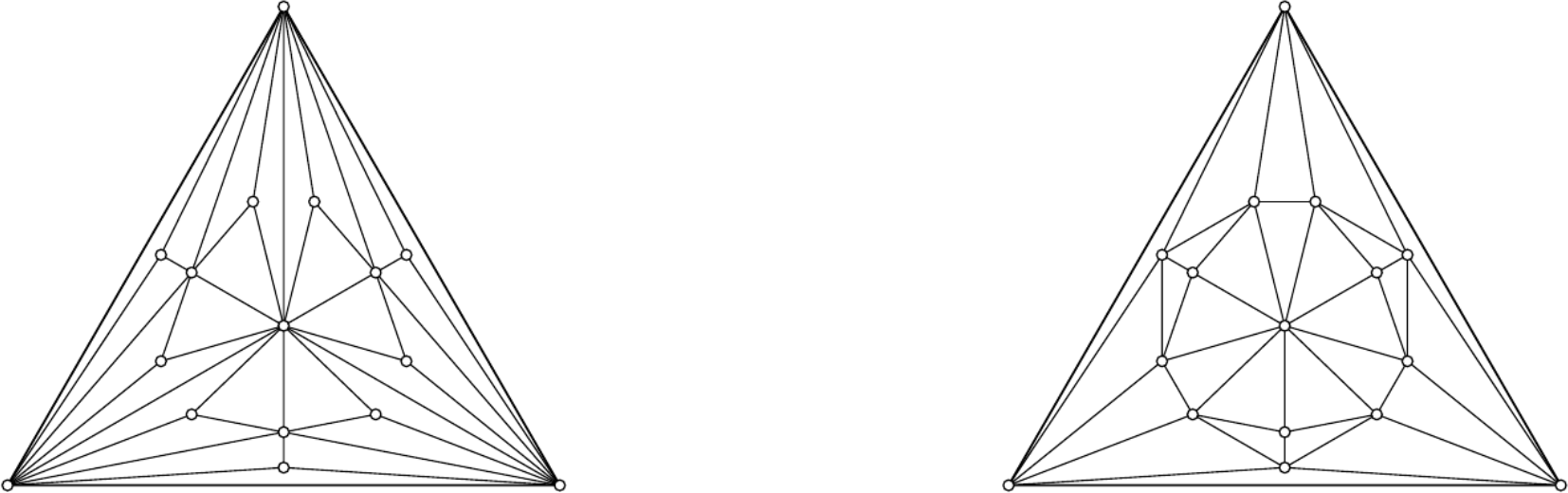}
\caption{(Left) three barycentric insertions of an equilateral triangle. (Right) the Delaunay triangulation of the same point set.}
\label{fig:delaunay}
\end{figure}

Next, our process stores a Delaunay triangulation of $S_0$. A \textit{Delaunay triangulation} is a triangulation of a point set which maximizes the minimum angle involved in a triangle (see  \cite{Triangulations} or \cite[Chapter~9]{computationalgeometry} for a definition). They are used ubiquitously in visualization algorithms due to this feature. In contrast, triangles with small angles often appear ``sharp'' or degenerate (see \Cref{fig:delaunay}). We rely on the \texttt{julia} package \texttt{DelaunayTriangulation.jl} for performing this standard computational geometry task efficiently \cite{Delaunayjl}.  

\subsection{Evaluation} Next, we evaluate the oracle $\mathcal O$ at each element of $S_0$. This is often the most computationally demanding step, motivating the main goal of the package: reduce the number of total evaluations required for a visualization by only evaluating at interesting regions of $D$. 

In practice, evaluating $f$ on many inputs may be best performed in batches. This is the case, for example, in \texttt{HomotopyContinuation.jl} \cite{HCjl} where it is much quicker to solve a parametrized polynomial system over many parameters in a single call than it is to run several single-solves.

\subsection{Completeness of Triangles}
\label{sub:completeness}
Initially, all triangles in the Delaunay triangulation are deemed \textit{incomplete}. The intuition behind this terminology is that incomplete triangles contain interesting behavior of $f$ and complete triangles do not. In practice, a triangle is complete if the values of $f$ on its vertices do not vary much. For discrete $f$, this means that the values are equal. For locally continuous $f$, it means that the values are within some numerical tolerance. 

Beyond these default methods for determining whether a triangle is deemed complete, an advanced user may provide their own \texttt{is\_complete} function.

\subsection{Main datatype: \texttt{TriangulationCache}}
The initial points $S_0$, their Delaunay triangulation, their evaluations via $\mathcal O$, and the completeness of each triangle are all stored in our main datatype of a \texttt{TriangulationCache}. A \texttt{TriangulationCache} can be easily initialized by calling it on any function of two variables. Optional arguments for its initialization include a visualization window (default: $[-1,1]^2$), a resolution (default: $100$), and several additional settings involved in the \textit{refinement} step discussed in the next section. 

\begin{example}
\label{ex:alphacurve}
The code below  produces an image (\Cref{fig:nodal}) of the Boolean function \[(x,y) \mapsto y^2 - x^2(x + 1) < 0.\]
\begin{lstlisting}
julia> nodal_cubic(x, y) = y^2 - x^2 * (x + 1) < 0 ? "negative" : "nonnegative"
julia> TC = TriangulationCache(nodal_cubic; xlims = [-2.5, 2.5],  resolution = 500)
julia> visualize(TC; edges=true)
\end{lstlisting} 
\begin{figure}[!htpb]
\includegraphics[scale=0.6]{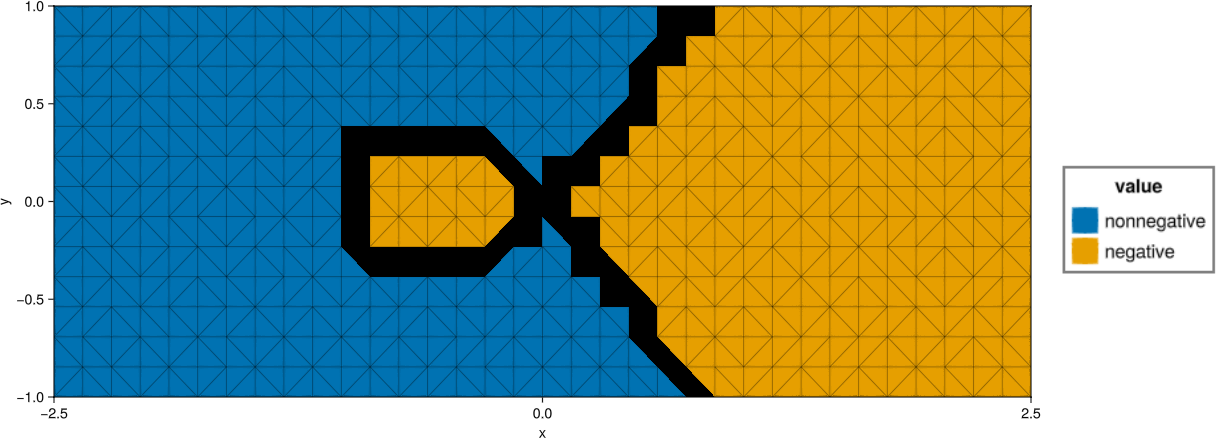}
\caption{Visualization of the function which determines if $f(x,y) = y^2-x^2(x+1)$ is negative. Called with resolution $N_0 = 500$.}
\label{fig:nodal}
\end{figure}

\end{example}

\begin{example}
\label{ex:petals}
Consider the function
\[(x,y) \mapsto (x^2 + y^2)^4 -
    (x^6 - 15x^4y^2 + 15x^2y^4 - y^6) < 0.\]
\Cref{fig:petal} illustrates the effect of using three different initial resolutions $N_0$: $500,5000,$ and $50000$.

\begin{figure}[!htpb]
\includegraphics[scale=0.34]{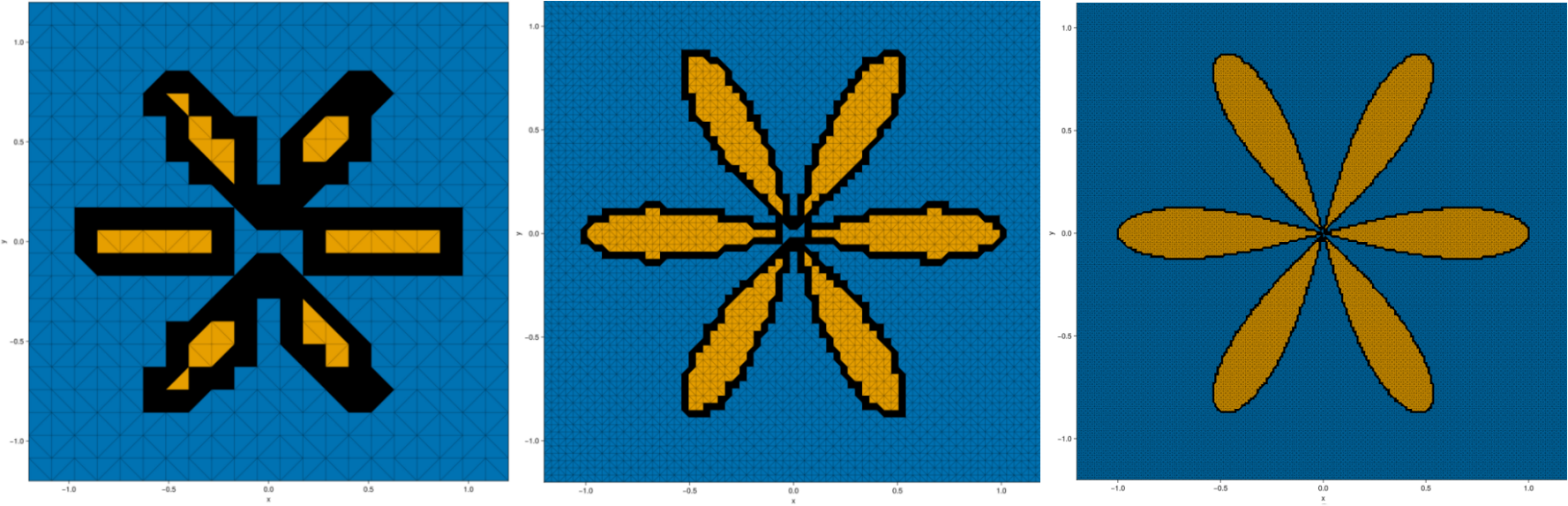}
\caption{Visualizations of the function which indicates whether a point lies inside or outside the curve cut out by  $(x^2 + y^2)^4 - (x^6 - 15x^4y^2 + 15x^2y^4 - y^6) = 0$. The resolutions of the visualizations are $500, 5000$, and $50000$, respectively.}
\label{fig:petal}

\end{figure}

\end{example}

\Cref{ex:petals} and \Cref{fig:petal} illustrate the crux of the problem with naive visualization via a uniform sample: increasing the resolution only marginally increases the visualization quality. Compare the rightmost image in \Cref{fig:petal} with the result of adaptive visualization,  \Cref{fig:example}, which only uses $5437$ evaluations.

\section{Refinement}
\label{sec:refinement}
\begin{figure}[!htpb]
$
\underset{\text{Select}}{
  \vcenter{\hbox{\includegraphics[scale=0.4]{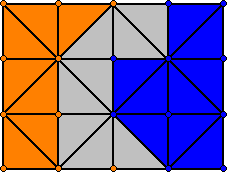}}}
}
\;\longrightarrow\;
\underset{\text{Insert}}{
  \vcenter{\hbox{\includegraphics[scale=0.4]{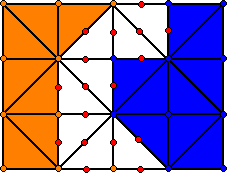}}}
}
\;\longrightarrow\;
\underset{\text{Triangulate}}{
  \vcenter{\hbox{\includegraphics[scale=0.4]{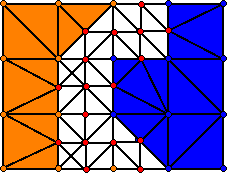}}}
}
\;\longrightarrow\;
\underset{\text{Evaluate}}{
  \vcenter{\hbox{\includegraphics[scale=0.4]{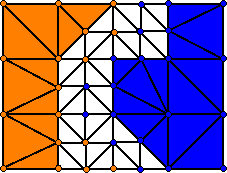}}}
}
\;\longrightarrow\;
\underset{\text{Completeness-check}}{
  \vcenter{\hbox{\includegraphics[scale=0.4]{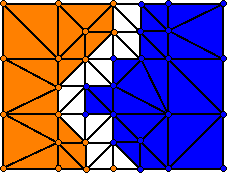}}}
}
$
    \caption{Illustration of the refinement steps. 
    }
    \label{fig:refinement_stage}
\end{figure}

The heart of our algorithm is the refinement step. Refinement involves inserting new points into each of the incomplete triangles, updating the Delaunay triangulation, and then repeating two of the steps already discussed in the previous section: ``evaluate'' and  ``completeness-check''.   The steps of the refinement stage \texttt{refine!} are illustrated in \Cref{fig:refinement_stage}. At the top level, \texttt{refine!} accepts a \texttt{budget} of a number of points that can be inserted. The process repeats until that budget runs out or no further refinement can occur.

\subsection{Selecting incomplete triangles} The first step of refinement is choosing  which incomplete triangles to refine. This is controlled by the window $D$ along with two keyword arguments: \texttt{min\_refinement\_area} and \texttt{budget}. The minimum refinement area controls one of the two filters which are always applied to the set of incomplete triangles. The first filter is that incomplete triangles outside of $D$ are ignored. This is relevant since our user-interface allows a user to move the viewing window. The second filter ignores incomplete triangles with area less than 
\[
\textrm{area}(D) \times \texttt{min\_refinement\_area}.
\] The default for \texttt{min\_refinement\_area} is $10^{-5}$, which, for reference, is about the proportion of a $1000 \times 1000$ image taken up by a $3 \times 3$ grid of pixels. 

Once the  candidate incomplete triangles have been established, they  are sorted by area and scheduled to be refined.

\subsection{Point insertion strategies} 
\begin{figure}[!htpb]
\includegraphics[scale=0.8]{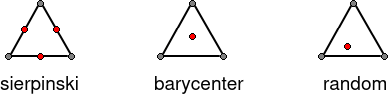}
\caption{Native local refinement strategies.}
\label{fig:strategies}
\end{figure}

Once a set of incomplete triangles has been chosen for refinement, points are sampled according to a \textit{local refinement strategy}. Our software offers three native choices:
\[
\begin{tabular}{rl}
\texttt{sierpinski:}& insert the edge midpoints of the inputted triangles, irredundantly.\\
 \texttt{barycenter:}& insert the barycenters of the inputted triangles.\\
 \texttt{random:}& insert random points in the interior of the inputted triangles.
\end{tabular}
\]
The last two of these increase the number of points in \texttt{TriangulationCache} by one for every incomplete triangle chosen for refinement. The method \texttt{sierpinski} increases the number of points by a factor of at most $3$. This depends on the configuration of the triangles chosen for refinement. These three strategies are illustrated in \Cref{fig:strategies}.

\subsection{Updating the Delaunay triangulation}
After the new points have been selected, they are inserted into the triangulation. \texttt{DelaunayTriangulation.jl} has functionality for inserting individual points into a Delaunay triangulation, triangulating, and then locally performing ``flips'' to make it Delaunay again. This locality is a well-known feature of Delaunay triangulations (see, e.g., \cite[Section 9.3]{computationalgeometry}) and is much more efficient than computing a brand new Delaunay triangulation.

\subsection{Evaluation and triangle completeness}

From here, the algorithm proceeds as it does at the end of the initialization step: $f$ is evaluated over each of the newly inserted points and each of the triangles in the updated \texttt{TriangulationCache} is checked for completeness.

\section{Functionality and User-Interface}
\label{sec:functionality}
We use this section to summarize the functionality of \texttt{AdaptiveVisualization.jl}.
The main exported functions of our package are $ \texttt{TriangulationCache}, \texttt{visualize},$ and $\texttt{refine!}$. We illustrate this functionality on the function from \Cref{ex:petals}. 

The \texttt{TriangulationCache} can be initialized via 
\begin{lstlisting}
julia> petal(x, y) = (x^2+y^2)^4-(x^6-15x^4*y^2+15x^2*y^4-y^6)<0 ? "inside" : "outside";
julia> I = [-1.2,1.2];
julia> TC = TriangulationCache(petal; xlims = I, ylims = I, resolution = 1000)
TriangulationCache with:
   961 function values
  1800 (306) triangles (incomplete)

julia> refine!(TC; budget = 7000); TC
Resolution used:589. Small triangles skipped:0
Resolution used:1237. Small triangles skipped:0
Resolution used:2534. Small triangles skipped:2507
Resolution used:110. Small triangles skipped:2627
Resolution used:6. Small triangles skipped:2640
Resolution used:0. Small triangles skipped:2640
Unable to spend refinement budget; no eligible refinement points remain.
TriangulationCache with:
   5437 function values
  10752 (2640) triangles (incomplete)
julia> visualize(TC)
\end{lstlisting}

\begin{figure}[!htpb]
\includegraphics[scale=0.3]{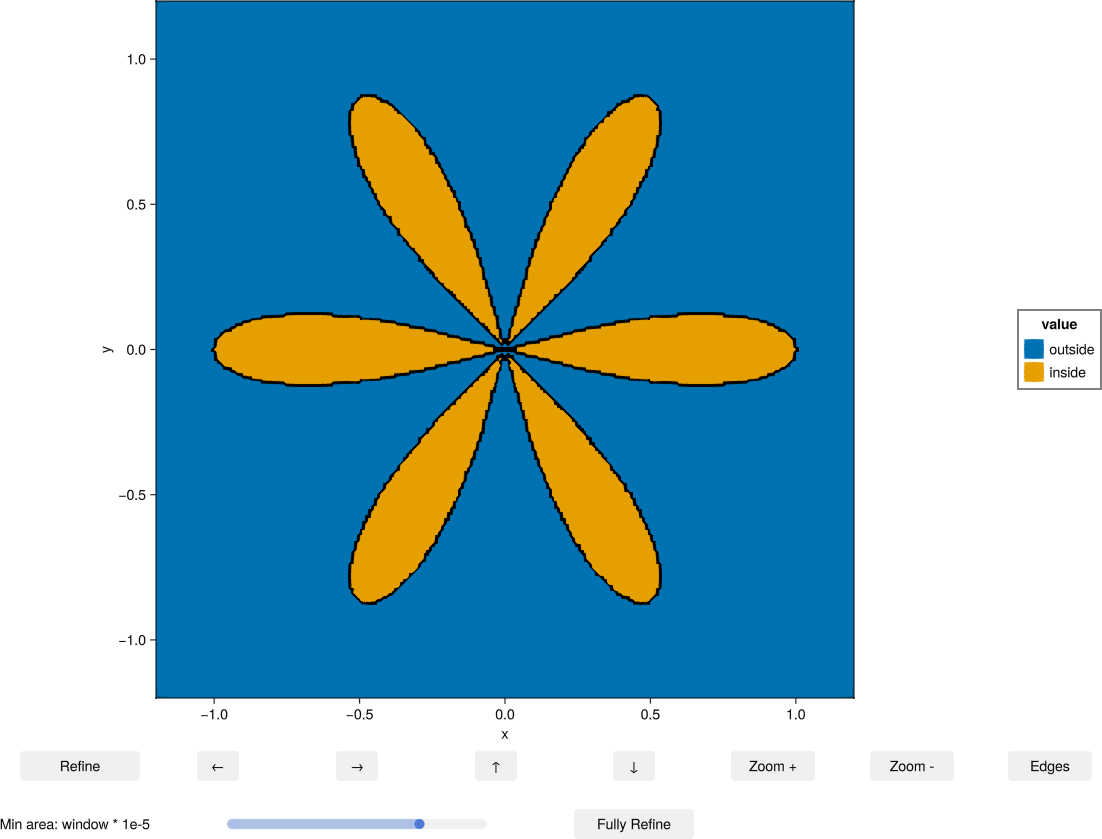}
\caption{Output of above code.}
\label{fig:example}
\end{figure}

The output Makie window is shown in \Cref{fig:example}.

Compare, for example, the quality of this image and that of the center image in \Cref{ex:petals} with approximately the same number of oracle calls ($5000$). 

\subsection{Makie Buttons} Notably, the Makie window shown in  \Cref{fig:example} has the buttons
\begin{center}
\texttt{Refine} \quad \quad $\leftarrow$ \quad \quad $\rightarrow$ \quad \quad $\uparrow$ \quad \quad $\downarrow$ \quad \quad \texttt{Zoom+} \quad \quad \texttt{Zoom-} \quad \quad \texttt{Edges} \quad \quad \texttt{Fully Refine}
\end{center} along with a slider for $\texttt{min\_refinement\_area}$.  \Cref{fig:petal} already illustrates how the edges of the triangulation appear when the \texttt{Edges} button is pressed. The remaining buttons do as one may expect: \texttt{Refine} calls \texttt{refine!(TC)} and replots, the buttons $\leftarrow, \rightarrow, \uparrow, \downarrow$ each translate the viewing window by $20\%$ in the corresponding direction, and \texttt{Zoom$\pm$} either zooms in or out by $20\%$. When the viewing window changes to contain regions of $\mathbb{R}^2$ not seen before, the code recognizes this and re-lays a lattice on the new window with resolution set to the optional argument \texttt{navigation\_initial\_resolution}. It then refines the new regions using a budget given by  the value \texttt{navigation\_refinement\_budget}. In particular, the software has the quirk that translating nearby a region repeatedly refines it. 
The $\texttt{min area}$ slider updates $\texttt{min\_area\_refinement}$. The \texttt{Fully Refine} button will iteratively call \texttt{refine!} until there are no incomplete triangles with sufficiently large area to refine. 

\begin{figure}[!htpb]
\includegraphics[scale=0.27]{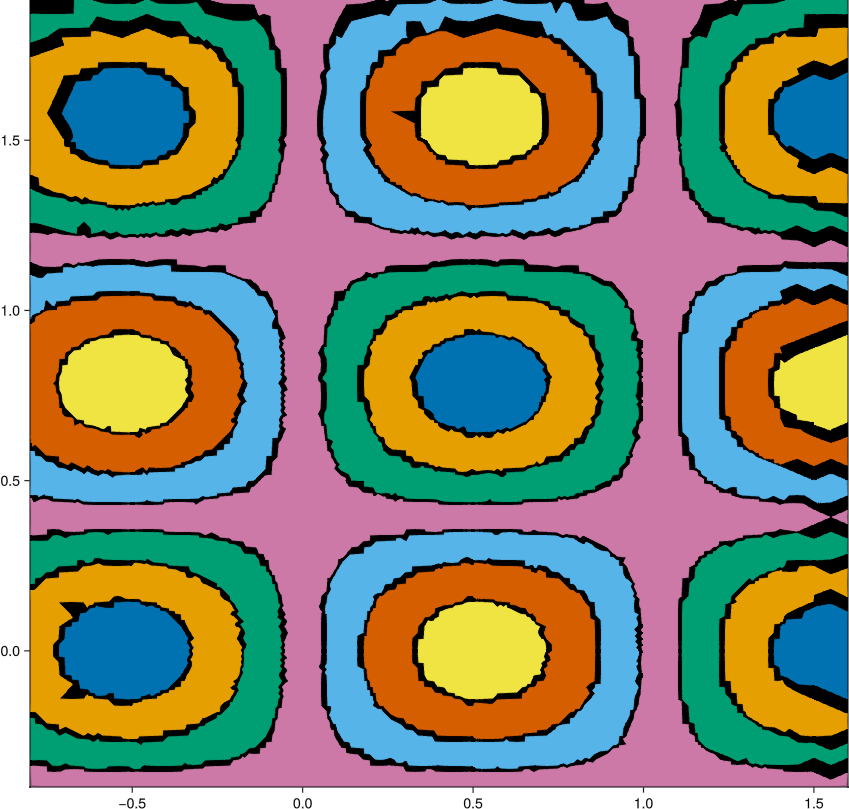}
\caption{Image of  $\text{round}(4 + 3\sin(3x)\cos(4y))$ after clicking $\rightarrow, \uparrow, \uparrow, \texttt{Zoom-}$ and \texttt{Refine}.}
\label{fig:puddles}
\end{figure}

The following code showcases the simplest way to use our software: pass an oracle to the \texttt{visualize} function. Doing so initializes a \texttt{TriangulationCache} with $\approx \frac 1 4 \texttt{total\_resolution}$-many points, and subsequently refines using the remaining budget of $\texttt{total\_resolution}$. 
\begin{lstlisting}
julia> f = (x, y) -> round(Int, 4 + 3sin(3x)cos(4y));
julia> (TC, puddleImage) = visualize(f)
Resolution used:598. Small triangles skipped:0
Resolution used:152. Small triangles skipped:0
(TriangulationCache with:
   975 function values
  1841 (978) triangles (incomplete), Scene(13 children, 0 plots))
\end{lstlisting}
The code above produces an image, which after pressing $\leftarrow, \uparrow, \uparrow, \texttt{Zoom-}, \texttt{Refine}$  produces \Cref{fig:puddles}.

\section{Applications to polynomial system solving}
\label{sec:enumerative_problems}
Our main motivation for the development of \texttt{AdaptiveVisualization.jl} was the use-case of visualizing parameter spaces of polynomial systems. Our software has a \texttt{julia} extension which activates  whenever both \texttt{AdaptiveVisualization} and \texttt{HomotopyContinuation} are loaded by the user. With this extension, the user can input  a polynomial system and explore its parameter space.

Fix a system of $n$ polynomials in $n$ complex variables and $k=2$ parameters,
\[
F(\x;\p) = \{f_1(\x;\p),\ldots,f_n(\x;\p)\} \subseteq \mathbb{C}[\p][\x] = \mathbb{C}[p_1,p_2][x_1,\ldots,x_n],
\]
with the property that $F(\x;\p^*)=\textbf{0}$ has finitely many complex solutions for generic $\p^*$. Then our software has hardcoded functionality to visualize the following three functions:

\begin{itemize}
\item[$N_{\mathbb{R}}$:] counts the number of real  solutions to $F(\textbf{x};p_1,p_2)=\textbf{0}$
\item[$N_{+}$:] counts the number of  positive real solutions to $F(\textbf{x};p_1,p_2)=\textbf{0}$
\item[$N_{\mathbb{R}, \checkmark}$:] the same as $N_{\mathbb{R}}$,  but certified using \texttt{HomotopyContinuation.jl}'s \cite{HCjl}  \texttt{certify} \cite{certify} command.\footnote{These are soft certificates since $N_{\mathbb{R}}$ is being evaluated at floating point input.} 
\end{itemize}

\begin{figure}[!htpb]
\includegraphics[scale=0.55]{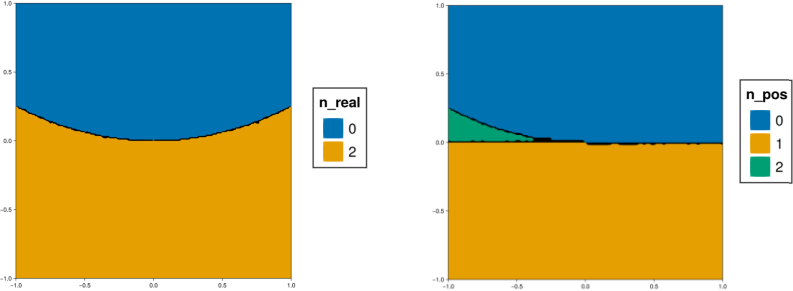}
\caption{Visualization of the number of real roots and the number of positive roots of a univariate quadratic polynomial.}
\label{fig:quadratic}
\end{figure}

 The simplest example is the number of real solutions to a monic quadratic univariate polynomial
\[
f(x;b,c) = x^2+bx+c.
\]
\begin{lstlisting}
julia> @var x; @var b,c;
julia> F = System([x^2+b*x+c],variables=[x],parameters=[b,c])
julia> TC1, fig1 = visualize(F; func=:real) # Default
julia> TC2, fig2 = visualize(F; func=:positive)
julia> TC3, fig3 = visualize(F; func=:certify_real) 
\end{lstlisting}
The pictures produced are shown in \Cref{fig:quadratic}. We remark that this computation is fundamentally different than how \cite{BreidingSturmfelsWang2025} compute similar images. 

The  locus $b^2-4c=0$, where the number of real solutions changes, emerges as the relief-image of the figures produced by our software. This is sometimes called the \textit{real discriminant}, which should not be confused with the real part of the complex discriminant (see \cite{MachineLearningDiscriminant} for an explanation). In settings with more solutions, it is known that chambers can be incredibly small (see, e.g., \cite[Page 4]{DickensteinEtAl2007}), making our methods particularly apt for such problems. 

One additional hard-coded function is the \textit{Dietmaier function} from \cite{Dietmaier}:
\[
\mathcal D: \mathbb{R}^2 \to \mathbb{R} \qquad \qquad 
(p_1,p_2) \mapsto\min_{\substack{\mathbf{x}\text{ a nonreal}\\\text{solution}}} (||\textrm{imag}(\textbf{x})||).
\]
It takes values in $\mathbb{R}_{\geq 0}$ and is the central function involved in the heuristic of \textit{Dietmaier optimization} (see 
 \cite{Dietmaier,DietmaierDavidDeepak}). 
Below, we visualize the Dietmaier function on the \textit{Kuramoto model} on three coupled oscillators (see \cite{Kuramoto1975,Mehta2015}). 

\begin{lstlisting}
julia> K = kuramoto_model(3); 
julia> visualize(K; func = :Dietmaier);
\end{lstlisting}

\begin{figure}[!htpb]
\includegraphics[scale=0.4]{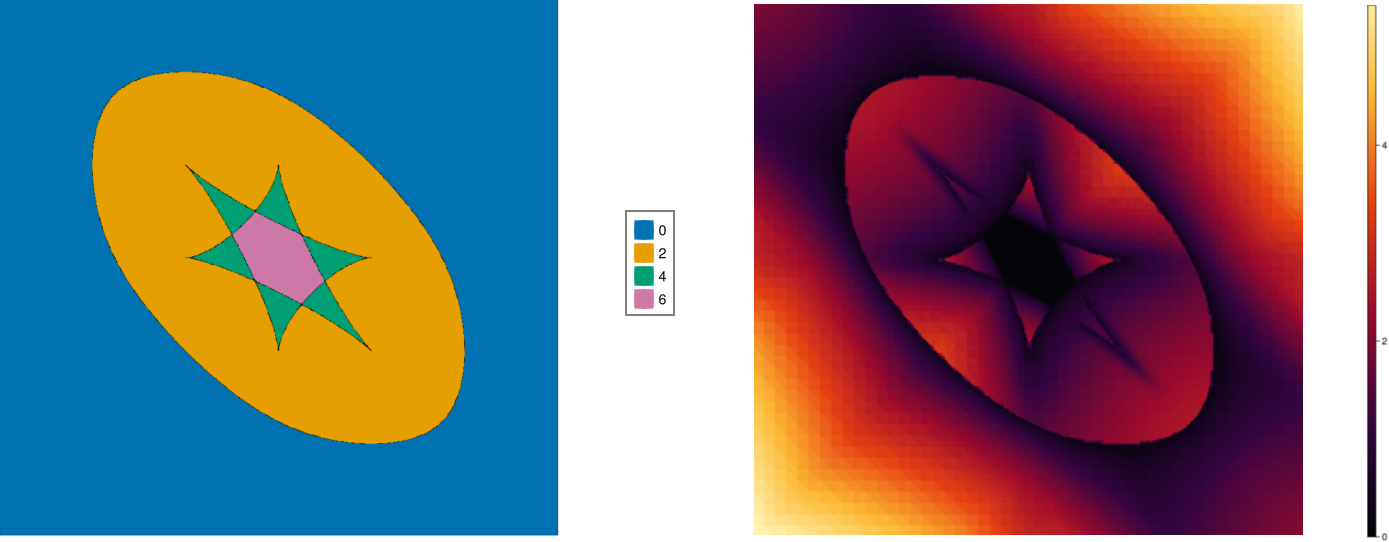}
\caption{(Left) The function which counts real solutions of the Kuramoto model for $n=3$. (Right) The Dietmaier function on the same model.} 
\end{figure}

\subsection{Kuramoto gallery} 
\begin{figure}[!htpb]
\includegraphics[scale=0.7]{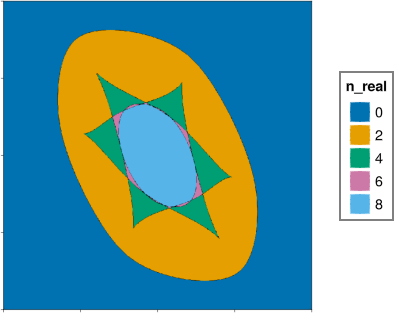}
\caption{A visualization of a random $2$-plane through $\textbf{0}$ in Kuramoto-4.}
\label{fig:kuramoto4}
\end{figure}

\begin{figure}[!htpb]
\includegraphics[scale=0.65]{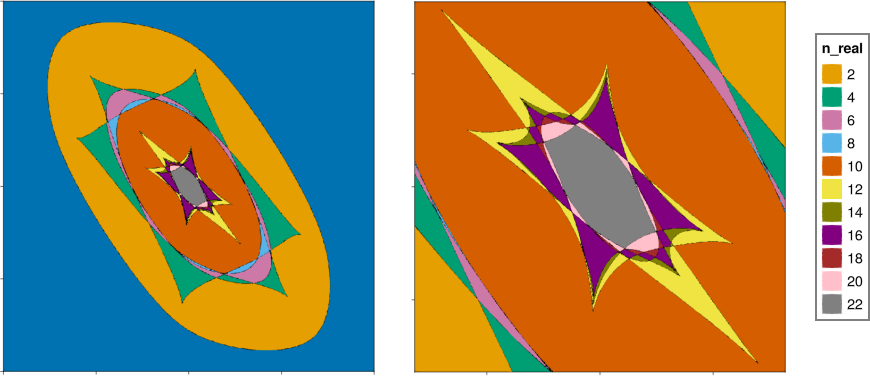}
\caption{A visualization of a random $2$-plane through $\textbf{0}$ in Kuramoto-5.}
\label{fig:kuramoto5}
\end{figure}

\begin{figure}[!htpb]
\includegraphics[scale=0.55]{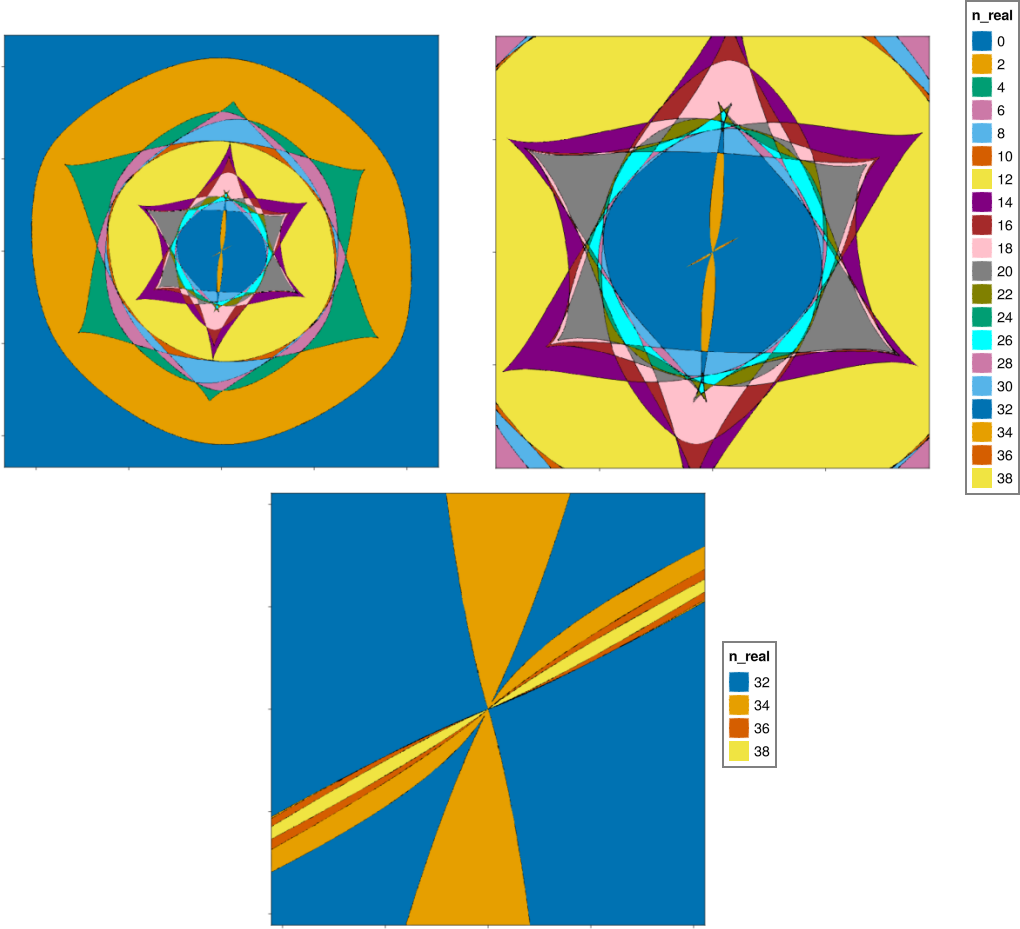}
\caption{A visualization of a random $2$-plane through $\textbf{0}$ in Kuramoto-6.}
\label{fig:kuramoto6}
\end{figure}

Kuramoto models for $n$ coupled oscillators (see \cite{Kuramoto1975,Kuramoto}) have beautiful real discriminants. We work with the algebraic formulation \begin{align*}
w_i-\frac{1}{n}\left(
s_i+\sum_{j=1}^{n-1}(s_i c_j-s_j c_i)
\right) &= 0,
&& i=1,\ldots,n-1,\\
c_i^2+s_i^2-1 &= 0,
&& i=1,\ldots,n-1.
\end{align*}
In \cite{MascarinTelen2026}, the authors relate these discriminants with Lissajous varieties. The authors of \cite{MachineLearningDiscriminant} also look at the real discriminant of the Kuramoto models, but they use intersection points with the discriminant and lines, where we do not ever directly access the discriminant algebraically or geometrically. 
In each image \autoref{fig:kuramoto4}-\autoref{fig:kuramoto6} we visualize a random $2$-plane through the origin of the parameter space of the Kuramoto model. The left-most plot is zoomed-out and the plots are zoomed-in as one moves to the right.

\begin{figure}[!htpb]
\includegraphics[scale=0.7]{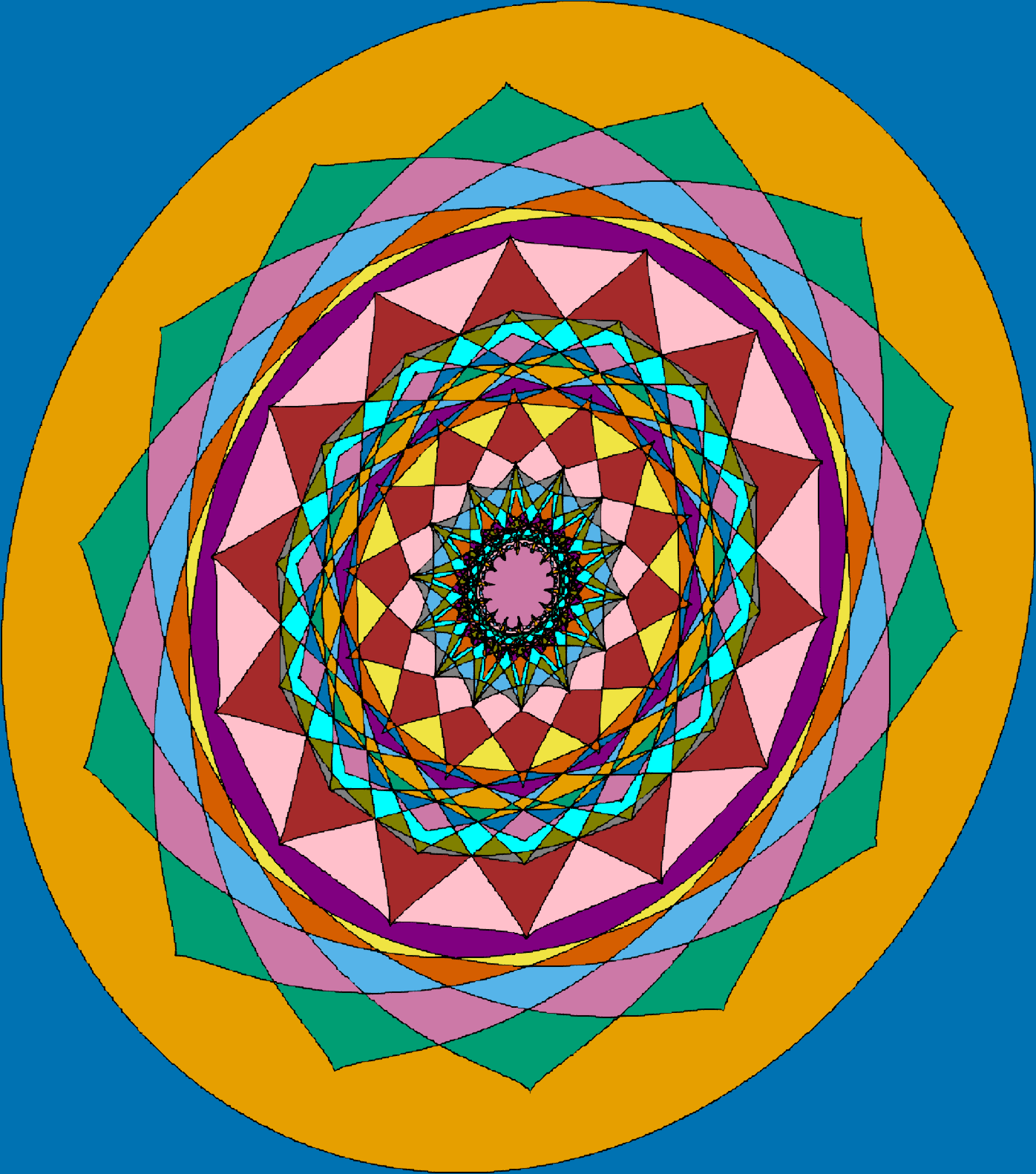}
\caption{A visualization of a  $2$-plane in the parameter space of Kuramoto-7, containing $\textbf{0}$, $ p = [0.70, 0.44, -0.16, -0.63, -0.63, -0.16]$ and $q = [-0.35, 0.33, 0.76, 0.62, 0.012, -0.60]$. This image has $14$-fold symmetry: $7$-fold symmetry by cyclically permuting the oscillators, and another $2$-fold symmetry given by negation. }
\label{fig:kuramoto7}
\end{figure}

We end with our favorite: a highly symmetric $2$-plane slice of the parameter space of Kuramoto-7.

\section{Acknowledgements and Funding Declaration} The authors are grateful to Charles Wampler for his comments on an initial version of the code and to Michael Burr for his comments on a later one. This project was supported by a USRA award and by NSERC Discovery Grant RGPIN-2023-03551.

\section*{AI Statement} AI tools were not used in the preparation of this manuscript. They were, however, used in improving the code. The original architecture and a prototype were produced by the authors alone. Later, with OpenAI's Astra and Sol models through Codex, the code was improved on the levels of user-interface and optimized. AI tools were used to simplify and organize the code as well as write tests and documentation. All code and documentation has been subsequently checked by the authors who take full responsibility for them.

\bibliographystyle{amsplain}
\bibliography{myrefs}

\end{document}